\documentclass[12pt]{article}

\usepackage[table]{xcolor}
\usepackage{amsmath, amsthm}
\usepackage{amsfonts,amssymb}
\usepackage{hyperref}
\usepackage{authblk}
\usepackage[all]{xy}
\usepackage{enumitem}
\usepackage{stmaryrd}
\usepackage[makeroom]{cancel}
\usepackage{bbm}
\usepackage{fancyvrb}

\usepackage[linesnumbered,boxed,vlined]{algorithm2e}

\usepackage{xstring} 

\usepackage{tikz}
\usetikzlibrary{intersections}

\newif\ifextmath
\extmathfalse

\newcommand{\N}{\mathbb{N}}
\newcommand{\Z}{\mathbb{Z}}
\newcommand{\Q}{\mathbb{Q}}

\newcommand{\C}{\mathbb{C}}

\newcommand{\F}{\mathbb{F}}
\renewcommand{\P}{\mathbb{P}}

\renewcommand{\Im}{\operatorname{Im}}

\newcommand{\lcm}{\operatorname{lcm}}
\newcommand{\codim}{\operatorname{codim}}
\newcommand{\eps}{\varepsilon}
\newcommand{\Fl}{\F_\ell}
\newcommand{\Flx}{\Fl^\times}
\newcommand{\GL}{\operatorname{GL}}
\newcommand{\SL}{\operatorname{SL}}
\newcommand{\PGL}{\operatorname{PGL}}
\newcommand{\PSL}{\operatorname{PSL}}
\newcommand{\GLFl}{\GL_2(\Fl)}
\newcommand{\SLFl}{\SL_2(\Fl)}
\newcommand{\PGLFl}{\PGL_2(\Fl)}
\newcommand{\PSLFl}{\PSL_2(\Fl)}
\newcommand{\Frob}{\operatorname{Frob}}
\newcommand{\disc}{\operatorname{disc}}
\newcommand{\Tr}{\operatorname{Tr}}

\newcommand{\Gal}{\operatorname{Gal}}

\newcommand{\mfurl}[1]{\StrSubstitute{#1}{.}{/}[\mfslash]}
\newcommand{\mfref}[1]{\mfurl{#1}\href{http://www.lmfdb.org/ModularForm/GL2/Q/holomorphic/\mfslash}{#1}}

\SetKw{True}{True}
\SetKw{False}{False}

\SetKw{FAIL}{FAIL}
\SetKw{Goto}{go to line}

\newtheorem{thm}[equation]{Theorem}
\newtheorem{lem}[equation]{Lemma}

\newtheorem{cor}[equation]{Corollary}
\theoremstyle{definition}
\newtheorem{de}[equation]{Definition}
\newtheorem{rk}[equation]{Remark}

\newtheorem{ex}[equation]{Example}

\makeatletter
\newcommand{\subjclass}[2][2020]{%
  \let\@oldtitle\@title%
  \gdef\@title{\@oldtitle\footnotetext{#1 \emph{Mathematics subject classification:} #2}}%
}
\newcommand{\keywords}[1]{%
  \let\@@oldtitle\@title%
  \gdef\@title{\@@oldtitle\footnotetext{\emph{Key words and phrases.} #1.}}%
}
\let\c@table\c@equation
\let\c@figure\c@equation
\makeatother

\makeatletter
\let\oldnl\nl
\newcommand{\nonl}{\renewcommand{\nl}{\let\nl\oldnl}}
\makeatother

\title{A method to prove that a modular Galois representation has large image}
\author{Nicolas Mascot\thanks{Trinity College Dublin, \href{mailto:mascotn@tcd.ie}{mascotn@tcd.ie}}}

\VerbatimFootnotes

\begin{document}

\maketitle

\begin{abstract}
Let~$\rho$ be a mod~$\ell$ Galois representation attached to a newform~$f$. Explicit methods such as~\cite{Samuele} are sometimes able to determine the image of~$\rho$, or even~\cite{MakMod} the number field cut out by~$\rho$, provided that~$\ell$ and the level~$N$ of~$f$ are small enough; however these methods are not amenable to the case where~$\ell$ or~$N$ are large. The purpose of this short note is to establish a sufficient condition for the image of~$\rho$ to be large and which remains easy to test for moderately large~$\ell$ and~$N$.
\end{abstract}

\renewcommand{\abstractname}{Acknowledgements}
\begin{abstract}
We express our thanks to B.~S.~Banwait for motivating us to write this short note by asking about the image of the representation discussed in Example~\ref{ex:Barinder} below.
\end{abstract}




\vspace{10mm}

Let~$f = q + \sum_{n \geqslant 2} a_n(f) q^n$ be a newform of weight~$k \in \N$, level~$N \in \N$, and nebentypus~$\eps : (\Z/N\Z)^\times \longrightarrow \C^\times$. Let~$\lambda$ be a prime of residual degree 1 of the field~$K_f$ spanned by the Hecke eigenvalues of~$f$, let~$\ell \in \N$ be the prime below~$\lambda$, and let~$\rho : \Gal(\overline \Q / \Q) \longrightarrow \GLFl$ be the Galois representation attached to~$f \bmod \lambda$. Since~$\rho$ also arises from a form of level prime to~$N$, we assume that~$\ell \nmid N$ from now on.

\begin{de}
Let~$p \in \N$ be a prime dividing~$N$. We say that~$f \bmod \lambda$ is \emph{old at~$p$} if~$\rho$ also arises from a form of level dividing~$N/p$.
\end{de}

Equivalently,~$f \bmod \lambda$ is old at~$p$ if it is congruent to a newform~$f'$ whose level divides~$N/p$; in particular,~$\eps$ must actually arise from a character mod~$N/p$. In fact, \cite[3.4, 7.2]{Edi} even show that this is equivalent to the existence of an integer~$i \in \N$, of a newform~$q + \sum_{n \geqslant 2} b_n q^n$ of weight~$k'$ satisfying~$2 \leqslant k' \leqslant \ell+2$, level dividing~$N/p$, and same nebentypus~$\eps$ as~$f$, and of a prime~$\lambda'$ above~$\ell$ of the field spanned by the~$b_n$ such that
\begin{equation}
a_n(f) \bmod \lambda = n^i b_n \bmod \lambda' \label{eqn:old}
\end{equation} for all~$n$ coprime to~$\ell N$. 
Therefore, given~$f$,~$\lambda$, and~$p$, it is possible (and usually not too difficult) to explicitly determine whether~$f \bmod \lambda$ is old at~$p$; see Examples~\ref{ex:Barinder} and~\ref{ex:PSL} below for concrete cases.

\begin{de}
We say that~$\rho$ has \emph{large} image if its image contains~$\SL_2(\F_\ell)$.
\end{de}

The purpose of this note is to establish an efficient one-way criterion for~$\rho$ to have large image. The key argument is established by the following Lemma.


\begin{lem}\label{lem:main}
Suppose that~$\rho$ is irreducible, and that there exists a prime~$p \mid N$ such that~$p^2 \nmid N$ and that the~$p$-part of~$\eps$ is trivial. Then either
\begin{itemize}
\item~$f \bmod \lambda$ is old at~$p$, or
\item~$\rho$ is ramified but not wildly ramified at~$p$, and the image by~$\rho$ of the inertia at~$p$ is cyclic of order~$\ell$.
\end{itemize}
\end{lem}

\begin{proof}
As $\rho$ is irreducible, Serre's conjecture (now a theorem thanks to~\cite{KW}) assigns to~$\rho$ a level~$N(\rho) \mid N$ and a weight~$k(\rho) \leqslant \ell^2-1$ such that~$\rho$ also arises from an newform~$f_\rho$ of level~$N(\rho)$ and weight~$k(\rho)$. While~$k(\rho)$ is defined in terms of the action of the inertia at~$\ell$,~$N(\rho)$ is defined in terms of the action of inertia at primes~$r \neq \ell$, viz
\[ N(\rho) = \prod_{r \neq \ell} r^{n_r}, \quad n_r = \codim V^{I_r} + \operatorname{Swan}_r, \]
where~$V \simeq \Fl^2$ is the representation space of~$\rho$,~$V^{I_r}$ is the subspace of~$V$ fixed by the inertia~$I_r$ at~$r$, and 
\[ \operatorname{Swan}_r = \sum_{i=1}^{+\infty} \frac1{[I_r:I_r^{(i)}]} \codim V^{I_r^{(i)}} \] is a non-negative integer defined in terms of the action on~$V$ of the higher inertia groups~$I_r^{(i)}$ at~$q$.

Therefore, if~$f \bmod \lambda$ is \emph{not} old at~$p$, then~$p \mid N(\rho)$. As~$N(\rho) \mid N$, it follows that
\[ 1 = n_p = \codim V^{I_p} + \operatorname{Swan}_p. \]
In particular,~$\codim V^{I_p} \leqslant 1$; but we cannot have~$\codim V^{I_p} = 0$ lest~$I_p$, and a fortiori the~$I_p^{(i)}$, act trivially on~$V$, which would also force~$\operatorname{Swan}_p = 0$. Therefore we must have~$\codim V^{I_p} =1$ and~$\operatorname{Swan}_p = 0$, meaning that up to conjugacy
\[ \rho_{\vert I_p} = \begin{pmatrix} 1 &  \xi \\ 0 & \chi \end{pmatrix} \]
for some morphism~$\chi : I_p \longmapsto \Flx$ and cocycle~$\xi : I_p \longmapsto \Fl$, and that~$\rho$ is ramified, but only moderately ramified, at~$p$. More specifically,~$\chi$ is the restriction to~$I_p$ of~$\det \rho$, and is therefore trivial: Indeed,~$\det \rho$ is the product of~$\eps$, which is trivial on~$I_p$ by our assumption that its~$p$-part is trivial, and of the~$(k-1)$-th power of the mod~$\ell$ cyclotomic character~$\chi_\ell$ (the one that expresses the Galois action on the~$\ell$-th roots of unity), which is also trivial on~$I_p$.

We deduce that actually
\[ \rho_{\vert I_p} = \begin{pmatrix} 1 &  \xi \\ 0 & 1 \end{pmatrix} \]
for some morphism~$\xi : I_p \longmapsto \Fl$, which cannot be trivial since~$\codim V^{I_p} \neq 0$.
\end{proof}

\begin{thm}\label{thm:main}
Suppose that there exists a prime~$p \mid N$ such that~$p^2 \nmid N$ and that the~$p$-part of~$\eps$ is trivial, and another prime~$r \nmid \ell N$ such that the polynomial~$x^2-a_r(f) x + r^{k-1} \eps(r) \bmod \lambda \in \Fl[x]$ is irreducible over~$\Fl$. Then either
\begin{itemize}
\item~$f \bmod \lambda$ is old at~$p$, or
\item~$\rho$ has large image.
\end{itemize}
\end{thm}

\begin{proof}
As the image of the Frobenius at~$r$ has irreducible characteristic polynomial~$x^2-a_r(f) x + r^{k-1} \eps(v) \bmod \lambda \in \Fl[x]$,~$\rho$ must be irreducible. Therefore, unless~$f \bmod \lambda$ is old at~$p$, Lemma \ref{lem:main} shows that that~$\rho(I_p)$ is cyclic of order~$\ell$, so~$\ell$ divides the order of~$\Im \rho$.
Proposition 15 of~\cite{SerEn} then shows that the image of~$\rho$ is either large or contained in a Borel subgroup of~$\GLFl$, but the latter case is incompatible with~$\rho$ being irreducible.
\end{proof}

\begin{rk}
Let $G \leqslant \GLFl$ be a subgroup containing $\SLFl$. The proportion of elements of~$G$ whose characteristic polynomial is irreducible only depends on~$\ell$; more specifically, it is~$\frac13$ if~$\ell=2$ and~$\frac12 \frac{\ell-1}{\ell+1}$ if~$\ell \geqslant 3$, which is always~$\geqslant \frac14$ and close to~$\frac12$ for large~$\ell$.

Therefore, if~$\rho$ indeed has large image, Cebotarev ensures that it will not be difficult to find a prime~$r \nmid \ell N$ such that the characteristic polynomial of~$\rho(\Frob_r)$ is irreducible. Conversely, it is reasonable to suspect $\rho$ is reducible if no such prime~$r$ is found after a few attempts; and if desired this suspicion may be rigorously confirmed by~\cite[Algorithm 7.2.4]{Samuele}.
\end{rk}

\begin{rk}
Suppose one wishes to search for pairs~$(f,\lambda)$ with~$\eps$ trivial and such that~$\rho$ has ``exotic'' image, meaning neither large nor contained in a Borel subgroup. One would presumably go through newforms~$f$ of increasing level~$N$, and not be interested in representations~$\rho$ already encountered in lower level. Theorem~\ref{thm:main} then shows that one can skip levels~$N$ having at least one non-repeated prime factor.
\end{rk}

\begin{rk}\label{rk:disc}
Suppose that the assumptions of Theorem \ref{thm:main} are satisfied, and that~$f \bmod \lambda$ is not old at~$p$, so that~$\rho$ has large image. Let~$k$ (resp.~$K$) be the number field of degree~$\ell+1$ (resp.~$\ell^2-1$) corresponding by~$\rho$ to the stabiliser of a point in~$\P^1(\Fl)$ (resp. a nonzero vector in~$\Fl^2$), and suppose~$N$ or~$\ell$ are too large for methods to determine~$k$ or~$K$ such as~\cite{MakMod} to apply, but that we hope to find a model for~$k$ or~$K$ in a database such as~\cite{LMFDB}.

Recall that whenever~$T(x) \in \Q[x]$ is irreducible of degree~$d$ and~$r \in \N$ is an at-most-tamely ramified prime in the root field~$F \simeq \Q[x]/T(x)$ of~$T(x)$, the exponent of~$r$ in the discriminant of~$F$ is~$d-\omega_r$, where~$\omega_r$ is the number of orbits of roots of~$T(x)$ under the action of inertia~$I_r$. Lemma \ref{lem:main} therefore also implies that the exponent of~$p$ in~$\disc k$ is~$\ell+1-2=\ell-1$, and it is~$\ell^2-1-2(\ell-1) = (\ell-1)^2$ in~$\disc K$. In particular, the primes above $p$ do not ramify in the extension $K/k$; this could already be seen in the proof of Lemma~\ref{lem:main}, which shows that $\rho(I_p)$ injects into $\PGLFl$.

Besides, when~$\ell \neq 2$, the fact that~$\rho$ is odd implies that the respective signatures of~$k$ and~$K$ are~$(2,\frac{\ell-1}2)$ and~$(\ell-1,\frac{\ell (\ell-1)}2)$. All this information, as well as the equivalence
\begin{align*}
a_r(f) = 0 \bmod \lambda \Longleftrightarrow & \Tr \rho(\Frob_r) = 0 \\
\Longleftrightarrow & \, \rho(\Frob_r) \text{ has order exactly 2 in $\PGLFl$} \\
\Longleftrightarrow & \, \text{$T_k(x)$ only has factors of degree 1 or 2 over~$\Fl$,} \\
& \text{with at least one factor of degree 2}  \\
\end{align*}

\vspace{-7mm}

\noindent whenever~$T_k(x) \in \Z[x]$ is an irreducible polynomial defining~$k$ and~$r \in \N$ is a prime not dividing~$\disc T_k$,
can help narrow down the search for~$k$ and~$K$ in a number field database.

Finally, if a candidate for~$k$ of~$K$ has been isolated, this identification can sometimes be proved rigorously by methods relying on Serre's conjecture, such as~\cite[4.2]{Companion}.
\end{rk}

Furthermore, our explicit knowledge of~$\det \rho$ makes it easy to explicitly determine the image of~$\rho$ when it is large:
\begin{cor}\label{cor:main}
Suppose that there exists a prime~$p \mid N$ such that~$p^2 \nmid N$ and that the~$p$-part of~$\eps$ is trivial, another prime~$r \nmid \ell N$ such that~$x^2-a_r(f) x + r^{k-1} \eps(r)$ is irreducible mod~$\lambda$, and that~$f \bmod \lambda$ is not old at~$p$. Let~$M = \lcm(\ell,N/p)$, and let
\[ \Delta = \{ x^{k-1} \eps(x) \bmod \lambda \ \vert \ x \in (\Z/M \Z)^\times \} \leqslant \Flx \]
and
\[ G =  \{ g \in \GLFl \ \vert \ \det g \in \Delta \} \leqslant \GLFl. \]
Then the image of~$\rho$ is~$G$.
In particular, if~$\ell \neq 2$, then the image of the projective representation attached to~$\rho$ is~$\PSLFl$ if~$\Delta$ is contained in the subgroups of squares of~$\Flx$, and~$\PGLFl$ else. (If~$\ell=2$, then it is~$\PGL_2(\F_2) = \PSL_2(\F_2)$, and $\Im \rho = \GL_2(\F_2) = \SL_2(\F_2)$.)
\end{cor}

\begin{proof}
This follows from Theorem \ref{thm:main} and for the fact that~$\det \rho = \chi_\ell^{k-1} \eps$, where~$\chi_\ell$ is the mod-$\ell$ cyclotomic character.
\end{proof}

\begin{ex}\label{ex:Barinder}
Let~$f=q-2q^4+(3+\sqrt2)q^5+O(q^6)$ be the newform of~\cite{LMFDB} label \mfref{9099.2.a.g}. Its level is~$N=9099 = 3^3\cdot p$ where~$p=337$ is prime, its weight is~$k=2$, its nebentypus is trivial, and its coefficient field is~$K_f=\Q(\sqrt2)$. Let~$\lambda$ be the prime~$(7,\sqrt2-3)$ of $K_f$, and let~$\rho$ be the corresponding representation.

Already for~$r=2$ we find that
\[ x^2-a_r(f) x + r = x^2+2 \]
is irreducible over~$\F_7$, which proves that $\rho$ is irreducible.

If~$f \bmod \lambda$ were old at~$p$, there would exist a newform~$f' = q + \sum_{n \geqslant 2} b_n q^n$ of level~$\Gamma_0(3^v)$,~$v \leqslant 3$ and weight~$2 \leqslant k' \leqslant 9$, and an integer~$i$ such that~\eqref{eqn:old} holds. As the \cite{LMFDB} informs us that~$a_r(f) \equiv 0 \bmod \lambda$ for~$r \in \{ 2,11,31,73 \}$, we run a computer search for such forms~$f'$ also satisfying that the~$b_r$ have norm divisible by~$7$ for all~$r \in \{ 2,11,31,73 \}$, which takes a couple of minutes and results in only possible form~$f'$. However, this~$f'$ also has~$b_5 \equiv 0 \bmod 7$ whereas~$a_5(f) \not \equiv 0 \bmod \lambda$, which proves that~$f$ is not old at~$p$.

By Corollary~\ref{cor:main}, we conclude that the representation attached to~$f \bmod \lambda$ is surjective. This answers a question of B.~S.~Banwait's, and which was our motivation for writing this short note.

The observations made in Remark \ref{rk:disc} show that in this example, the field~$k$ is a number field of signature~$(2,3)$, unramified away form~$\{3,7,337\}$ and certainly ramified at~$7$ (by~$\det \rho$) and tamely ramified 337, the exponent of~$337$ in~$\disc k$ being~$6$; and the Galois group of the Galois closure of~$k/\Q$ is~$\PGL_2(\F_7)$. This information is more than enough to determine that this field is not present in the~\cite{LMFDB} as of May 2022.
\end{ex}

\begin{rk}
On the other hand,~$f$ is old at $p$ modulo the conjugate prime~$\lambda' = (7,\sqrt2+3)$, as $f \bmod \lambda'$ is congruent to the newform of weight~$2$, level~$27$, and trivial nebentypus. In fact, this form of weight~$2$ happens to have CM, and accordingly the image of the projective representation attached to~$f \bmod \lambda'$ is dihedral according to~\cite[Table 2]{Barinder}.
\end{rk}

\begin{ex}\label{ex:PSL}
This time, let~$f = q + \sum_{n \geqslant 2} a_n(f) q^n$ be the newform of~\cite{LMFDB} label~\mfref{71.3.b.a}. It has prime level~$N=71$, weight~$k=3$, quadratic nebentypus~$\eps = \left( \frac{\cdot}{71} \right )$, and coefficient field~$K_f = \Q(\beta)$ where~$\beta^4+108\beta^2-40\beta+2825=0$. Let us take for instance~$\lambda = (41,\beta-11)$, a prime of~$K_f$ of degree~$1$ above~$\ell=41$, and let~$\rho$ be the corresponding mod~$\ell$ Galois representation.

First of all, we check that for~$r=2$, the characteristic polynomial
\[ x^2-a_r(f) x + r^{k-1} \eps(r) \bmod \lambda = x^2-16x+4 \]
is irreducible over~$\Fl$, so~$\rho$ is irreducible.

Let us thus apply Theorem \ref{thm:main} with~$p=N$. If~$f \bmod \lambda$ were old at~$p$, there would exist a newform~$f' = q + \sum_{n \geqslant 2} b_n q^n$ of level~$1$ and weight at most~$\ell+2=43$ and an integer~$i$ such that \eqref{eqn:old} holds. However, a computer search reveals in less than one second that~$a_{101}(f) = 0 \bmod \lambda$ whereas all these forms~$f'$ have that the norm of~$b_{101}$ is not a multiple of~$\ell$. Therefore~$f \bmod \lambda$ is not old at~$p$, so~$\rho$ has large image.

In fact, as~$\det \rho = \chi_\ell^2 \eps$ and as the values~$\pm1$ of~$\eps$ are squares mod~$\lambda$, Corollary~\ref{cor:main} shows that the image of~$\rho$ is the subgroup of~$\GL_2(\F_{41})$ formed of matrices whose determinant is a square mod~$41$, and the corresponding projective representation has image~$\PSL_2(\F_{41})$.

In particular, in this example, the Galois closure of~$k/\Q$ has Galois group the simple group~$\PSL_2(\F_{41})$, and ramifies exactly at~$\ell$ and at~$p$, tamely at~$p$, and probably wildly at~$\ell$.
\end{ex}

\end{document}